\documentclass[letterpaper,10 pt, conference]{ieeeconf}
\IEEEoverridecommandlockouts
\usepackage{amsmath,amssymb,amsfonts}
\usepackage{algorithm}
\usepackage[noend]{algpseudocode}
\usepackage{graphicx}
\usepackage{textcomp}
\usepackage[dvipsnames,table]{xcolor}
\usepackage{comment}
\usepackage{booktabs}
\usepackage{url}
\usepackage{makecell}
\usepackage{subcaption}
\usepackage[para]{footmisc}

\let\labelindent\relax
\usepackage[inline]{enumitem}
\usepackage{siunitx}
\usepackage[printonlyused,withpage]{acronym}

\makeatletter
\let\NAT@parse\undefined
\makeatother
\usepackage[hidelinks]{hyperref}
\usepackage[capitalise,nameinlink,sort]{cleveref}

\def\BibTeX{{\rm B\kern-.05em{\sc i\kern-.025em b}\kern-.08em
    T\kern-.1667em\lower.7ex\hbox{E}\kern-.125emX}}

\newcommand{\argmin}{\operatornamewithlimits{argmin}}
\newcommand{\R}{\mathbb{R}}
\newcommand{\X}{\mathcal{X}}

\graphicspath{{figs/}}
\title{ \textbf{Data-Driven Extrusion Force Control Tuning for 3D Printing}}

\author{Xavier Guidetti$^{1,2}$, Ankita Mukne$^{1}$, Marvin Rueppel$^{2}$, Yannick Nagel$^{3}$, Efe C. Balta$^{1,2}$, and John Lygeros$^{1}$%
\thanks{Research supported by Innosuisse (project \textnumero 102.617 IP-ENG) and by the Swiss National Science Foundation under NCCR Automation (grant \textnumero 180545).  $^1$Automatic Control Laboratory (IfA), ETH Z\"urich, 8092 Z\"urich, Switzerland. $^2$Control and Automation Group, Inspire AG, 8005 Z\"urich, Switzerland. $^3$NematX AG, 8093 Z\"urich, Switzerland. \newline
Corresponding author: E. C. Balta, \texttt{ \href{mailto:efe.balta@inspire.ch}{efe.balta@inspire.ch}}.}
}

\begin{document}

\maketitle

\begin{abstract}
The quality of 3D prints often varies due to different conditions inherent to each print, such as filament type, print speed, and nozzle size. Closed-loop process control methods improve the accuracy and repeatability of 3D prints. However, optimal tuning of controllers for given process parameters and design geometry is often a challenge with manually tuned controllers resulting in inconsistent and suboptimal results. This work employs Bayesian optimization to identify the optimal controller parameters. Additionally, we explore transfer learning in the context of 3D printing by leveraging prior information from past trials. By integrating optimized extrusion force control and transfer learning, we provide a novel framework for closed-loop 3D printing and propose an automated calibration routine that produces high-quality prints for a desired combination of print settings, material, and shape.
\end{abstract}

\acrodef{fcp}[FCP]{Force Controlled Printing}

\section{Introduction}
Additive Manufacturing (AM), often known as 3D printing, enables manufacturing of complex geometries using a layerwise bottom-up process methodology.
A popular method for AM in practice is Fused Filament Fabrication (FFF) with polymers. 
FFF uses a heated extruder system to melt polymer material and print it in the form of ellipsoidal beads in a layerwise fashion.
The recent development of high-performance technical polymers allows FFF printed parts to have high mechanical strength and stiffness, expanding their use in a wide range of application areas~\cite{guidetti2023stress}.
Additionally, fiber-reinforced polymer printing is commonly employed in various advanced manufacturing applications~\cite{goh2019recent}.

A critical challenge in FFF is ensuring the reliability and repeatability of the process under changing process conditions. 
On most current applications, the challenge is partially caused by a lack of in-situ sensors, control-oriented models, and closed-loop control strategies~\cite{balta2021layer}.
A possible approach to improve performance is through modeling printed bead outputs as a function of input parameters~\cite{balta2022numerical,aksoy2020control}.
Similarly, extrusion dynamics models can be used for feedforward input optimization to improve dimensional accuracy~\cite{wu2023modeling}.
Even though open-loop model-based approaches are powerful and beneficial for improving process performance, model inaccuracies and run time disturbances remain a challenge. 
As a result, in-situ measurement of run-time performance for monitoring and control is necessary for high-performance AM in practice.


In-situ monitoring for FFF processes has received attention in the recent literature~\cite{rao2015online}.
Applications of health management~\cite{kim2018development}, quality monitoring~\cite{chhetri2019quilt}, and anomaly detection~\cite{balta2019digital} have been successfully developed for FFF using in-situ measurements.
Meanwhile, process control received comparatively less attention. 
An effective strategy for closed-loop control is using measurements of the printed part between layers. Layer-to-layer measurements of FFF prints have been used for accurate process modeling~\cite{balta2021layer} and parameter optimization~\cite{guidetti2023data}.
However, layer-to-layer updates are ineffective for rejecting fast-acting run-time disturbances to the extrusion flow in the layer, which has a great influence on the resulting printed material shape~\cite{balta2022numerical,serdeczny2018experimental,serdeczny2019numerical}.
Recent developments consider in-layer measurements of the material flow within the extruder~\cite{coogan2019line}.
Such measurements provide valuable run-time feedback and can be used for closed-loop process control.
In this work, we characterize the extrusion flow through feedback force measurements and develop a novel in-layer closed-loop control framework for FFF.

Controlling the extrusion force enables an effective way to build controllers that ensure desired flow characteristics.
However, developing and tuning such controllers is challenging due to the complex characteristics of the extrusion process coupled with the printing dynamics. 
Manually tuned controllers cannot guarantee optimal results. 
Additionally, the nonlinear process dynamics are difficult to model for real-time model-based control.
In this work, we employ Bayesian Optimization, a data-driven optimization method, to find optimal control parameters for \acl{fcp} \cite{guidetti2024force} for FFF and present effective Transfer Learning strategies for knowledge transfer between optimal tunings. The main contributions of this work are:
\begin{itemize}[leftmargin=*]
    \item The development of a preliminary \acl{fcp} framework for FFF;
    \item A continuous Bayesian Optimization method for controller parameter tuning during a print process;
    \item The use of Transfer Learning on Bayesian Optimization to accelerate convergence rates.
\end{itemize}
We experimentally demonstrate our results on a custom-built setup to showcase the effectiveness. 
Our work provides a baseline for future developments in force-controlled FFF.

\section{Background} \label{sec2}

\subsection{Force Controlled Fused Filament Fabrication}

\begin{figure}[htbp]
\centerline{\includegraphics[trim={0 0.3cm 0 0.2cm},clip]{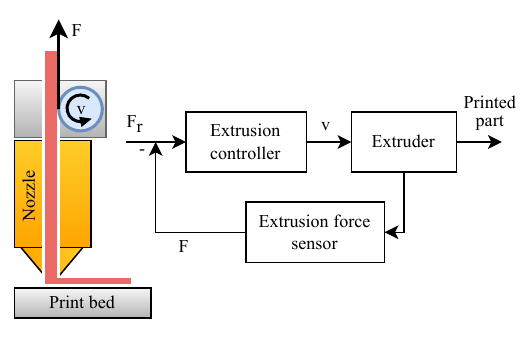}}
\caption{Simplified diagram of the extrusion process and feedback loop for \acs{fcp} for FFF}
\label{fig:freebody}
\end{figure}

In a conventional FFF extruder, a thermoplastic filament is pushed by a driving gear into a heated nozzle. The extruder is moved over the build plane where the material exiting the nozzle is deposited to form a part. The volume of extruded material is pre-computed based on the filament size, the nozzle diameter, and the extruder motion~\cite{aksoy2020control}. 
In \acf{fcp} for FFF, a sensor is installed in the extruder assembly to measure the force applied onto the filament while extruding it through the nozzle during printing. 
Unlike the conventional approach, where the filament driver speed is predefined, in \acs{fcp} the driver speed is continuously adapted by a controller to maintain the measured driving force at a desired and constant reference value (see \cref{fig:freebody}). The printing performance of this technique strongly depends on the controller's performance. PID controllers~\cite{visioli2006practical} can be used for this application, but only a well-tuned PID controller will properly regulate the extrusion process while dealing with external disturbances.

\subsection{Bayesian Optimization}

BO is a data-driven iterative optimization method, used to optimize expensive-to-evaluate black-box functions that can be sampled only in a point-wise fashion. It is best suited for optimization over continuous domains and tolerates stochastic noise in function evaluations~\cite{frazier2018tutorial, garnett2023bayesian}. BO has been successfully used in a variety of applications such as machine learning \cite{snoek2012practical}, robotics \cite{berkenkamp2023bayesian}, control \cite{khosravi2021performance}, manufacturing \cite{guidetti2022advanced}, and more; \cite{guidetti2023data} uses BO to tune the printing process parameters of FFF. 

Consider a problem of optimizing an unknown objective $f: \mathcal{X} \rightarrow \R$ over a set of inputs $\mathcal{X}\subset\R^d$:
\begin{equation}
    \min_{x \in \X} f(x).
    \label{eq:bboopt}
\end{equation}
In each iteration we observe a noisy evaluation $\textstyle{y_t= f(x_t) + \varphi_t}$ for the chosen input $x_t \in \mathcal{X}$,
where $\varphi_t \sim \mathcal{N}(0,\rho^2)$ is zero-mean noise, independent across different time steps $t$.
We collect the past input-output data in the set $\textstyle{D_t=\{(x_i,y_i)\}_{i=1}^t}$.
Using the collected data in $D_t$, BO uses a Gaussian Process (GP)~\cite{williams2006gaussian} to model the unknown function and an acquisition function to choose subsequent samples.

GPs are commonly used in non-parametric modeling as they provide a distribution over functions from limited data. Let $\mu_t(\cdot)$ and $\sigma^2_t(\cdot)$, denote the posterior mean and variance, respectively, given the previous measurements $\mathbf{y}_t = [y_1, \dots, y_t]^\top$ and kernel $k(\cdot,\cdot)$. These are computed by
\begin{align}
    \mu_t(x) &= k_t(x)^T(K_t + \rho^2I)^{-1} \mathbf{y}_t \label{eq:GPmean},\\
     \sigma^2_t(x) &=  k(x,x) - k_t(x)^\top(K_t + \rho^2I)^{-1}k_t(x) \label{eq:GPvar},
\end{align} 
where $(K_t)_{i,j} = k(x_i, x_j)$, $\ k_t(\cdot)^T = [k(x_1, \cdot), .., k(x_t, \cdot) ]^T$.


The acquisition function 
$\mathrm{acq_t}: \X \rightarrow \R$ 
is used for understanding how informative an input is for finding the minimizer of~\eqref{eq:bboopt} given the GP model of $f$. 
A good acquisition function effectively trades off exploration of new uncertain inputs with the exploitation of observed data that leads to better iterates with lower cost. 
\cref{alg:bo} provides the standard BO loop. 
\begin{algorithm}[ht]
\caption{Bayesian Optimization Loop}
 \label{alg:bo}
    \begin{algorithmic}[1] 
       
        \State \textbf{Initialize}: Prior $f \sim GP(0, k)$, $\bar{x} = 0$
        \For{$t$ = 1, \dots, $T$}
            \State $x_t \leftarrow \argmin\{ \mathrm{acq_t}(x)~|~x \in \X\}$
            \State Observe $y_t \gets f(x_t) +\varphi_t$
            \State Update the GP posterior~\eqref{eq:GPmean}-\eqref{eq:GPvar} using $y_t$
            \If{$f(x_t) < f(\bar{x})$} $\bar{x} \gets x_t$
            \EndIf
        \EndFor
        \State \textbf{Output}: $\bar{x}$
    \end{algorithmic}
\end{algorithm}

A stopping condition (such as a maximum number of iterations) is commonly used to interrupt the BO. One of the key strengths of BO is its efficiency in finding the optimal solution with a relatively small number of evaluations. 

\subsection{Transfer Learning}
Transfer learning (TL) is a machine learning technique in which any knowledge learned on one task is reused to improve performance on a secondary related task. For example, TL is very popular in deep learning, where pre-trained models are used to warm start the training process of a specialized model \cite{tan2018survey}. 
TL improves or accelerates the understanding of
the current task by relating it to other tasks performed at different periods through a related source domain~\cite{hosna2022transfer}. 
TL can be integrated with Bayesian Optimization, by modifying the surrogate model, allowing the outcomes of one optimization task to inform the optimization process of a related task. 
TL has been used in FFF by \cite{zhu2024surface}, where it was combined with a convolutional neural network for surface feature prediction. In our framework, we use TL to accelerate the convergence of BO on new tasks.

\section{Proposed Framework}\label{sec3}

\subsection{Controller Structure}

We implement a modified PID controller to regulate the extrusion force. The proportional and integral gains correspond to a common PID structure. However, we modify the derivative part to include two derivative terms. These two terms are designed to compute numerical derivatives of the force measurements at different time resolutions: 
\begin{enumerate}[label=\alph*),leftmargin=*]
    \item Fast (more responsive) derivative accounting for the faster machine motion system with a small time constant;
    \item Slow (more damped) derivative accounting for the slower filament extrusion system, having a larger time constant.
\end{enumerate}
By processing the measurements at the two different time constants, we are able to control the printing process more effectively. 
Ultimately, our PID controller has four tunable gains: $K_p$ (the proportional gain), $K_i$ (the integral gain), $K_d$ (the fast derivative), and $K_{dd}$ (the slow derivative).

\subsection{Continuous Bayesian Optimization} \label{sec:method:bo}

\begin{figure}[htbp]
\centerline{\includegraphics[trim={0 0.5cm 0 0.7cm},clip]{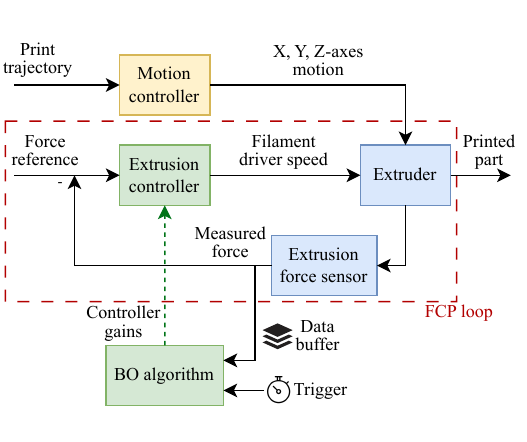}}
\caption{Flowchart of the continuous Bayesian optimization method}
\label{fig:contBO}
\end{figure}

We implement BO (see~\cref{alg:bo}) to find the optimal gains of the PID controller regulating the extrusion force. 
Specifically, to learn the optimal controller for a specific geometry, filament, printing temperature, and nozzle we conduct the BO routine \emph{while} printing an exact copy of the desired part. Unlike most BO applications -- where each iteration corresponds to one experiment, such as printing an entire part -- we define each iteration as a \SI{10}{s} long controller deployment executed while printing a part. As manufacturing a part by FFF is a time-consuming process, we can conduct a large number of BO iterations during one print. We call this approach \emph{continuous BO}. 

In practice, the process, which we illustrate in \cref{fig:contBO}, works as follows: 
\begin{enumerate}
    \item The printing instructions for a desired geometry (G-Code) are executed, and the printer head moves along a predefined trajectory;
    \item Extrusion control is initialized by setting the desired force reference and a random initial controller;
    \item The software records \SI{10}{s} of printing data with the given controller; the BO routine evaluates the printing performance and immediately deploys a new controller; this point is repeated until completion of the printing instructions or until user interruption.
\end{enumerate}
This procedure continuously employs~\cref{alg:bo} on an FFF system in a single print.
The system keeps printing with no interruption during the time needed in between iterations to evaluate the controller performance, update the surrogate models, run the acquisition function, and deploy the next controller. 
This means that discovering the optimal controller for a given geometry is achieved with one single print requiring the same printing time as the final part.
By the end of a single \SI{10}{\minute} experiment, 60 controllers have been tested, and the controller with the best performance will be used to print the desired part in a subsequent print. 

Next, we provide the specific formulation of BO used in our framework. The feedback control scheme aims to minimize the error between the measured extrusion force and a desired reference force. Given a reference force $F_r$ and a sequence of consecutive force measurements $\{F_i\}_{i=1}^n$, we compute the force root-mean-square error (RMSE) as
\begin{equation}
\mathrm{RMSE}(\{F_i\}_{i=1}^n,F_r) = \sqrt{\frac{1}{n}\sum_{i=1}^n (F_i - F_r)^2} \, .
\end{equation}
We evaluate the performance of a controller by computing the RMSE of the data recorded while the controller is deployed. Since the ideal controller has $\mathrm{RMSE} = 0$, we aim to minimize the RMSE. Considering a set of controller gains $x = \{K_p,K_i,K_d,K_{dd}\}$, and the force recorded while deploying this controller $\{F_i(x)\}_{i=1}^n$, 
we define the BO black-box objective function as
\begin{equation}
f(x) = \mathrm{RMSE}(\{F_i(x)\}_{i=1}^n,F_r) \, . \label{eq:objBO}
\end{equation}
We select expected improvement (EI) as an acquisition function~\cite{zhan2020expected}, a popular strategy for information collection which can be computed analytically for a GP surrogate model. EI is defined as
\begin{equation}
\mathrm{EI}(x) = \mathbb{E}[\max\{f(\bar{x}) -f(x) +\xi , 0\}] \, ,   
\end{equation}
where $\bar{x}$ represents the best controller discovered so far and $\xi$ is a tuning parameter influencing the selections made by the acquisition function. Larger values of $\xi$ lead to more exploration, while smaller values of $\xi$ produce more exploitation. By experimenting with the FCP setup, we find a suitable trade-off by setting $\xi = \SI{1e-3}{}$. We bound the inputs by $0 \leq K_p \leq 50$, $0 \leq K_i \leq 50$, $0 \leq K_d \leq 100$, and $0 \leq K_{dd} \leq 100$. The collected data is modeled using a GP with a Matern kernel having $\nu = 2.5$. The hyperparameters of the GP kernel are selected by fitting the GP on the available data and maximizing the log marginal likelihood. 

\subsection{Transfer Learning with Bayesian Optimization}

We consider the case where we want to print the same geometry with different reference forces. Naturally, each reference force corresponds to a different optimal controller. The goal of TL is to accelerate and improve the discovery of the optimal controller corresponding to a new reference force by exploiting the data that was previously collected with other reference forces; we refer to the different reference forces as ``different tasks".  

To utilize TL in conjunction with BO, we made minor modifications to the formulation of the BO discussed in \cref{sec:method:bo}. First, we extend the input vector $x$ by appending the reference force, leading to $x = \{K_p,K_i,K_d,K_{dd},F_r\}$. Then, to produce a more uniform data set for better task learning, the objective function is normalized using the force reference to reduce inter-task variations: 
\begin{equation}
f(x) = \frac{\mathrm{RMSE}(\{F_i(x)\}_{i=1}^n,F_r)}{F_r} \, .
\end{equation}
The acquisition function (EI) remains unchanged.

\section{Experimental Study}\label{sec4}

\subsection{Hardware Setup and Printing Conditions}

Our \acs{fcp} setup comprises a printer connected to three different sub-systems. The displacement of the extruder in space is controlled by a programmable logic controller (PLC) using Beckhoff\footnote{\url{https://www.beckhoff.com/}} software (corresponding to the yellow block in \cref{fig:contBO}). This unit operates the motors moving the printer head along the three Cartesian axes. The sequence of motions that are executed depends on a printing instructions file 
loaded into the Beckhoff software. 
The second unit has a force sensor and a Raspberry Pi board actuating the filament extruder stepper motor (corresponding to the blue blocks in \cref{fig:contBO}). 
The last unit is a computer running ROS2~\cite{ros2}
(corresponding to the green blocks in \cref{fig:contBO}), which closes the feedback loop for \acs{fcp} and runs the PID controller in run-time. 
The force sensor embedded in the extruder streams real-time data to the ROS2 unit. 
Based on the measurements and the control and/or optimization algorithms, ROS2 streams new extrusion values or parameters to the Raspberry Pi, which in turn updates the extruder drive speed.

All experiments were conducted by printing the \emph{tower} shells shown in \cref{fig:printed_shells}, hollow cubes with no top or bottom, made of a single and continuous contour line. 
This geometry was selected as the motion and printing conditions repeat identically throughout the print, making the learning and optimization tasks more effective.
Prints were made at an axis speed of \SI{100}{mm/s}, temperature of \SI{300}{\celsius}, layer height of \SI{50}{\micro\meter}, using a \SI{1.75}{mm} nozzle, and while extruding a liquid crystal polymer produced by NematX AG\footnote{\url{https://nematx.com/}}.

\subsection{Results}

\subsubsection{Single-Task Bayesian Optimization}
\label{sec:res:bo}

\begin{figure}[htbp]
    \centering
    \begin{subfigure}[t]{0.48\columnwidth}
    \includegraphics[width=\columnwidth]{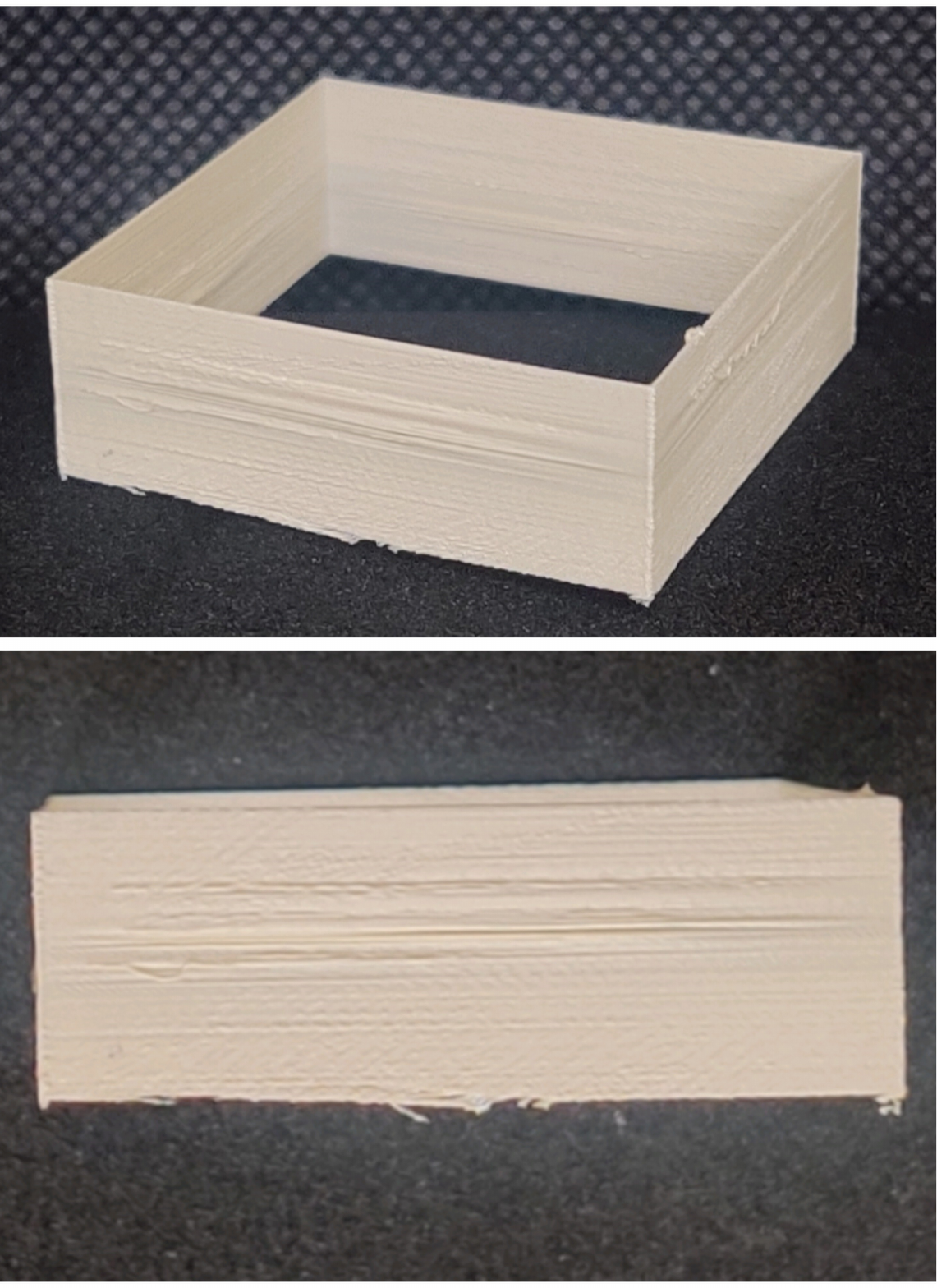}
    \caption{Before BO}\label{subfig:beforeBO}
    \end{subfigure}
    \begin{subfigure}[t]{0.48\columnwidth}
    \includegraphics[width=\columnwidth]{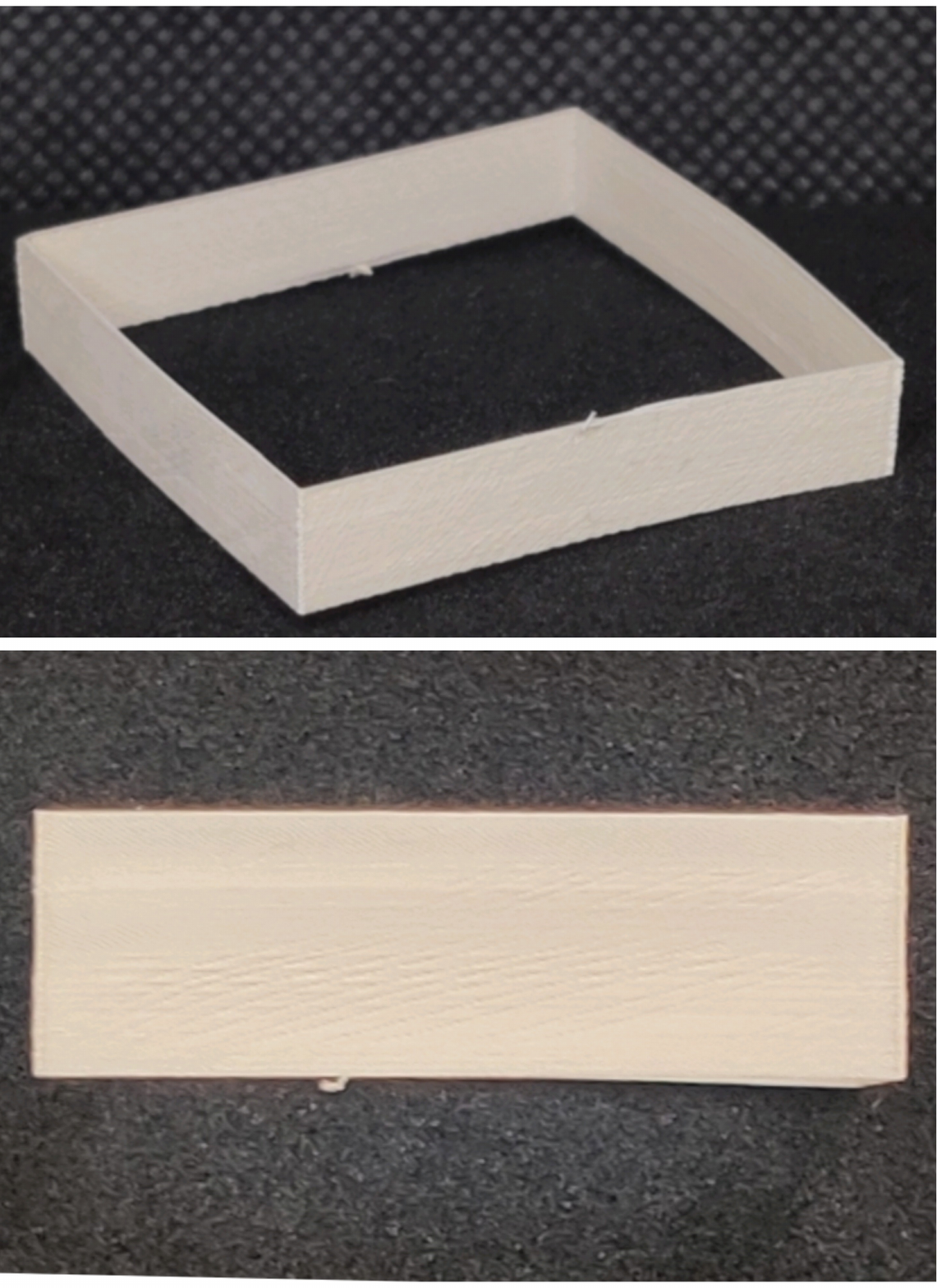}
    \caption{After BO}  \label{subfig:afterBO}
    \end{subfigure}
    \caption{Printed \emph{tower} shells on the experimental setup with force feedback extrusion control. 
    Both parts are continuously printed until user interruption, with force reference \SI{0.3}{N}.}
    \label{fig:printed_shells}
\end{figure}

The \emph{tower} was printed while running single-task continuous BO in four different experiments, each with a different reference force. The results are summarized in \cref{tab:singleBOres}. In all cases, the performance of a controller manually tuned by a practitioner was used to initialize the optimization. Then continuous BO was executed, and the performance of the best-found controller was recorded. The results show a significant reduction in RMSE for all cases, which corresponds to an increase in print quality. This observation can be confirmed qualitatively by analyzing \cref{fig:printed_shells}. The figure shows two towers printed under identical conditions at a force reference of \SI{0.3}{N} before and after controller optimization. In~\cref{subfig:beforeBO}, it is clearly visible how the side walls of the part have a rough surface. This is because machine's abrupt velocity and direction changes at the corners of the part disturb the extrusion process. As the manually tuned controller is unable to properly reject these disturbances, numerous defects appear in the \emph{tower} walls. After BO, however, the extrusion controller performance is drastically increased. In \cref{subfig:afterBO}, the walls show practically no surface defects, since the controller regulates the extrusion process more effectively.

\begin{table}[htbp]
\centering
\caption{Comparison of controller performance before and after Bayesian optimization}
\label{tab:singleBOres}
\renewcommand{\arraystretch}{1.3}
\begin{tabular}[h]{@{} c c c @{}}
\toprule
Reference [N] & RMSE before BO [N] & RMSE after BO [N]\\
\midrule
0.1 & \SI{4.029e-3}{} & \SI{3.477e-3}{} \\
0.2 & 0.289 & 0.038 \\
0.3 & 0.284 & 0.047 \\
0.4 & 0.200 & 0.011 \\
\bottomrule
\end{tabular}
\end{table}

\begin{figure*}[htbp]
    \centering
    \begin{subfigure}[t]{0.32\textwidth}
    \includegraphics{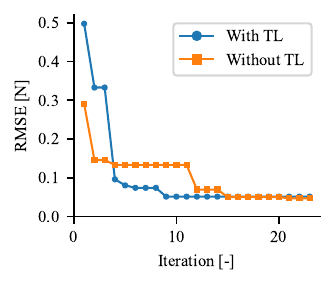}
    \caption{$F_r = \SI{0.3}{N}$}\label{fig:tl_4_3}
    \end{subfigure}
    \begin{subfigure}[t]{0.32\textwidth}
    \includegraphics{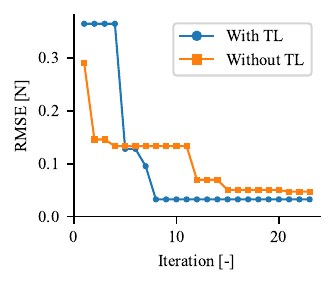}
    \caption{$F_r = \SI{0.3}{N}$}  \label{fig:tl_2_3}
    \end{subfigure}
    \begin{subfigure}[t]{0.32\textwidth}
    \includegraphics{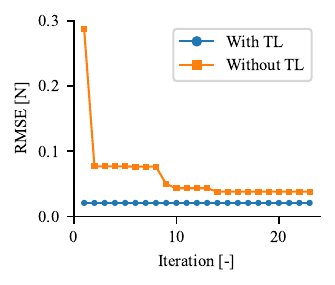}
    \caption{$F_r = \SI{0.2}{N}$}  \label{fig:tl_4_3_2}
    \end{subfigure}
    \caption{Comparison of convergence plots from controller optimization for: (a) force reference of \SI{0.3}{N} with TL using data from the optimization of \SI{0.4}{N}, (b) force reference of \SI{0.3}{N} with TL using data from the optimization of \SI{0.2}{N}, (c) force reference of \SI{0.2}{N} with TL using data from the optimization of both \SI{0.3}{N} and \SI{0.4}{N}.}
    \label{fig:comp_bo}
\end{figure*}

\subsubsection{Transfer Learning} \label{sec:res:tl}

To study the effect of TL on BO, we conducted a set of experiments where controllers were tuned for a new force reference while having access to data belonging to the optimization conducted previously for another force reference. In particular, we were interested in observing whether the use of TL would produce a controller with better performance and reduce the number of iterations required to optimize the process. In the three experiments we conducted, we
\begin{enumerate*}[label=\emph{\alph*})]
    \item used optimization data from a force reference of \SI{0.4}{N} to warm-start the optimization for a force reference of \SI{0.3}{N}, \label{exp_a}
    \item used optimization data from a force reference of \SI{0.2}{N} to warm-start the optimization for a force reference of \SI{0.3}{N}, and \label{exp_b}
    \item used optimization data from force references of \SI{0.4}{N} \emph{and} \SI{0.3}{N} to warm-start the optimization for a force reference of \SI{0.2}{N}. \label{exp_c}
\end{enumerate*}

Figures \ref{fig:tl_4_3} and \ref{fig:tl_2_3} show the results for experiments \ref{exp_a} and \ref{exp_b}, comparing the convergence of the optimization using TL for a target task of reference \SI{0.3}{N} with the results produced without TL. In both cases, the initial iterations produce worse performance, as commonly observed in TL. The optimization process trusts excessively the structure learned in the training task, which does not necessarily generalize to the target task. However, a few samples collected in the target task are enough for the TL approach to capitalize on the inter-task similarities, leading to faster convergence. In addition, the performance of the optimal controller is either similar to the one produced without TL (\cref{fig:tl_4_3}), or better (\cref{fig:tl_2_3}). The main metrics are compared quantitatively in \cref{tab:tl_comparison}.

\begin{table}[htbp]
\centering
\caption{Comparison of cost (RMSE) and required iterations between different controller tuning strategies, for a target force reference of \SI{0.3}{N}}
\label{tab:tl_comparison}
\renewcommand{\arraystretch}{1.3}
\begin{tabular}[h]{@{}l c c @{}}
\toprule
& Best RMSE [N] & \makecell[c]{Iterations \\ to convergence}\\
\midrule
Without TL & 0.047  & 21  \\
With TL (from \SI{0.4}{N})& 0.052  & 9  \\
With TL (from \SI{0.2}{N})& 0.033  & 8  \\
\bottomrule
\end{tabular}
\end{table}

In \cref{fig:tl_4_3_2} we show the results from experiment \ref{exp_c}, where data from optimization for references \SI{0.4}{N} and \SI{0.3}{N} are used to warm-start BO for reference \SI{0.2}{N}. Note that the target task is outside the range of tasks used in the training, requiring the TL to extrapolate. We compare the convergence of the optimizations conducted with and without TL. The results clearly show that the model learned from two training tasks predicts exceptionally well the target task: the optimal controller is discovered at the first iteration and it significantly outperforms the controller discovered by conventional BO. Quantitative metrics are reported in \cref{tab:double_tl_comp}.

\begin{table}[htbp]
\centering
\caption{Comparison of cost (RMSE) and required iterations between different controller tuning strategies, for a target force reference of \SI{0.2}{N}}
\label{tab:double_tl_comp}
\renewcommand{\arraystretch}{1.3}
\begin{tabular}[h]{@{}l c c @{}}
\toprule
& Best RMSE [N] & \makecell[c]{Iterations \\ to convergence}\\
\midrule
Without TL & 0.038  & 14\vspace{3pt}  \\
\makecell[l]{With TL \\ (from \SI{0.4}{N} and \SI{0.3}{N})} & 0.020  & 1 \\
\bottomrule
\end{tabular}
\end{table}

\begin{figure}[htbp]
\centerline{\includegraphics{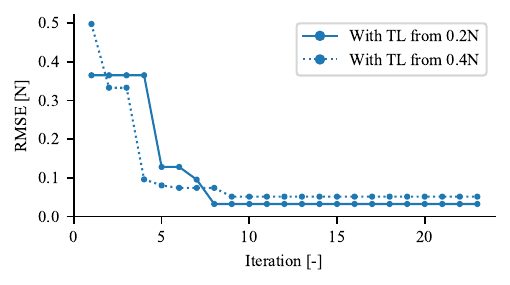}}
\caption{Comparison of convergence plots from controller optimization for a force reference of \SI{0.3}{N}}
\label{fig:tl_2&4}
\end{figure}

\subsection{Discussion} 

Our results in \cref{sec:res:bo} demonstrate how for the force references of \SI{0.1}{N}, \SI{0.2}{N}, \SI{0.3}{N}, and \SI{0.4}{N} the RMSE of manually tuned controllers was reduced by \SI{13.7}{\percent}, \SI{86.9}{\percent}, \SI{83.5}{\percent}, and \SI{94.5}{\percent} respectively by using standard BO. This translates directly into better-quality parts. It appears from the results that BO is a powerful method for automated controller tuning in \acs{fcp}. 
In the experiment with reference \SI{0.1}{N}, BO appears produce subpar results, as the RMSE improvement is not as significant as in the other tasks. This is due to the fact that the manually tuned controller used for benchmarking was already producing excellent results, since regulating the extrusion force at low references appears to be relatively simple. Nonetheless, BO was able to efficiently discover a better controller.

The results of \cref{sec:res:tl} show how TL in conjunction with BO can further improve performance, significantly accelerating the optimization and leading to better controllers. The improvements are summarized in \cref{tab:improvements}. In only one case, TL produced a marginally worse controller. This was however compensated by a drastic reduction in the number of iterations required for convergence. In \cref{fig:tl_2&4} we compare the two cases where TL was conducted with the same target task (\SI{0.3}{N} reference) but with two different training tasks (\SI{0.2}{N} and \SI{0.4}{N}). Using data from the \SI{0.2}{N} reference, TL converges marginally faster and to a better controller. We believe this is due to the characteristics of the controllers found in the training data set. In particular, the data set generated with a lower reference (\SI{0.2}{N}) contains controllers with lower proportional and integral gains (data not shown). This seems to skew the controller exploration in favor of less aggressive control gains, making the optimization more efficient. Conversely, the optimization using the higher reference (\SI{0.4}{N}) data set was observed to test numerous unstable controllers, which affected the optimization negatively.

\begin{table}[htbp]
\centering
\caption{Change in cost (RMSE) and iterations produced by different transfer learning strategies, when compared to controller tuning without transfer learning}
\label{tab:improvements}
\renewcommand{\arraystretch}{1.3}
\begin{tabular}[h]{@{}l r r @{}}
\toprule
& Best RMSE &\makecell[c]{Iterations \\ to convergence}\\
\midrule
\SI{-0.3}{N} from \SI{-0.4}{N}  & \SI{10.6}{\percent}  & \SI{-57.1}{\percent}   \\
\SI{-0.3}{N} from \SI{-0.2}{N}   & \SI{-29.8}{\percent}  & \SI{-61.9}{\percent}\vspace{1pt}  \\
\makecell[c]{\SI{-0.2}{N} from \\ \SI{-0.4}{N} and \SI{-0.3}{N}}  & \SI{-47.4}{\percent} & \SI{-92.8}{\percent}   \\
\bottomrule
\end{tabular}
\end{table}

\section{Conclusion and Future Work} \label{sec5}

We introduced a framework for force controlled FFF and demonstrated experimentally how continuous BO can be used to efficiently tune extrusion controllers that optimize the printing performance. The method requires a single print and reliably produces high-performance controllers. We have also utilized TL in conjunction with BO and shown that the optimization is significantly faster and leads to better controllers when transferring information from past experiments. Future work focuses on characterizing the mechanical properties of parts printed with optimized \acs{fcp} and on extending the TL framework to changes in part geometry.

\bibliographystyle{IEEEtran}
\bibliography{ref.bib}

\end{document}